\documentclass{amsproc}

\usepackage{amsmath}
\usepackage{amscd}
\usepackage{amssymb} 

\newcommand{\cal}{\mathcal}
\newcommand{\bk}{{\bf k}}

\newcommand{\bC}{{\Bbb C}}

\newcommand{\bQ}{{\Bbb Q}}
\newcommand{\bR}{{\Bbb R}}

\newcommand{\cA}{{\cal A}}
\newcommand{\cB}{{\cal B}}

\newcommand{\cG}{{\cal G}}
\newcommand{\cH}{{\cal H}}
\newcommand{\cL}{{\cal L}}

\newcommand{\fd}{{\frak d}}
\newcommand{\fm}{{\frak m}}

\DeclareMathOperator{\Img}{Im}

\DeclareMathOperator{\Ker}{Ker}

\theoremstyle{remark}

\theoremstyle{definition}
 \newtheorem{example}{Example}[section]

\begin{document}
\title
{DGBV Algebras and Mirror Symmetry}
\author{Huai-Dong Cao, Jian Zhou}
\address{Department of Mathematics\\
Texas A\&M University\\
College Station, TX 77843}
\email{cao@math.tamu.edu, zhou@math.tamu.edu}
\begin{abstract}
We describe some recent development on the theory of formal Frobenius 
manifolds via a construction from differential Gerstenhaber-Batalin-Vilkovisk 
(DGBV) algebras and formulate a version of mirror symmetry conjecture:
the extended deformation problems of the complex structure and 
the Poisson structure are described by two DGBV algebras;
mirror symmetry is interpreted in term of the invariance of the formal
Frobenius manifold structures under quasi-isomorphism.
\end{abstract}
\maketitle

\section{Introduction}

According to Getzler \cite{Get}, topological string theory is conformal 
field theorist's algebraic topology. Indeed,
ideas from cohomology theory has also been used extensively by physicists, 
especially in the two widely used quantization schemes:
the BRST formalism and the BV formalism.
On the other hand,
infinite algebra structures and the notion of operad
originally developed in homotopy theory,
have also shown up in many places in string theory.
A well-known connection between the cohomology theory and the homotopy theory
is provided by the rational homotopy theory.
One naturally speculates that such a connection should have its counterpart in
string theory.
Here we report some recent work on mirror symmetry
which reflects this connection.

The mirror symmetry is one of the mysteries in string theory.
For its history and backgrounds,
see Yau \cite{Yau1} and Greene-Yau \cite{Gre-Yau}.
In physicist's language,
given a Calabi-Yau manifold,
one can define two kinds of superconformal field theories on it:
the A-type theories and the B-type theories.
The mirror symmetry conjecture says for certain Calabi-Yau manifolds $M$,
there exist Calabi-Yau manifolds $\widehat{M}$,
such that a certain A-type theory on $M$
can be identified with a certain B-type theory on $\widehat{M}$.
There are three issues involved in this conjecture:
\begin{itemize}
\item[(a)] the construction of the mirror manifolds;
\item[(b)] the mathematical formulations of the A and B-type  theories;
\item[(c)] the identifications of the relevant theories.
\end{itemize}
Most of researches so far has focused on (a) and (b).
We present here an approach which deals with (b) and (c) in the same framework.

According to physicists,
an A-type theory should be sensitive to the deformation 
of the K\"{a}hler structure,
but is independent of the complex structure,
while a B-type theory should be sensitive 
to the deformation of the complex structure,
but is independent of the K\"{a}hler structure.
So it is reasonable to speculate that such theories are related to 
the deformation theories of the K\"{a}hler 
structure and the complex structure respectively.
It will be clear from below that it is more natural to consider 
the deformation of the Poisson structure defined by the K\"{a}hler structure.

In deformation theory,
the following principle due to Deligne
is well-known (see e.g. Goldman-Millson \cite{Gol-Mil}):
``In characteristic zero a deformation problem
is controlled by a differential graded Lie algebra with the property that
quasi-isomorphic differential graded Lie algebras give the same deformation
theory."
In this principle, 
one finds many objects in the rational homotopy theory.
The notion of quasi-isomorphism was used in Sullivan's minimal model
theory approach to rational homotopy theory \cite{Sul}.
As suggested by the BV quantization scheme,
the consideration of the extended deformation problem is necessary.
As a modification of Deligne's principle,
one notices that such problems are 
usually governed by a differential Gerstenhaber (or Poisson) algebra.
For example,
the extended deformation problems of the complex structure 
and the Poisson structure of a Calabi-Yau manifold
are controlled by 
two differential Gerstenhaber-Batalin-Vilkovisky (DGBV) algebras 
in the title.
Such algebras are combinations of differential graded algebras (DGA's) 
and differential graded Lie algebras (DGLA's).
Chen \cite{Che} developed a theory for DGA's 
to compute the (co)homology of loop spaces. Hain \cite{Hai} generalized it 
to DGLA's and showed that Chen's theory is an alternative to Sullivan's theory.
Our approach to (b) uses a construction of formal Frobenius manifolds
from DGBV algebras in which the formal power series connection on
a DGLA is used. 
We have defined in \cite{Cao-Zho5} a notion of quasi-isomorphisms of DGBV 
algebras, and have shown that formal Frobenius manifolds obtained from 
quasi-isomorphic DGBV algebras can be identified with each other.
This is our approach to (c).

\section{Frobenius algebras and formal Frobenius manifolds}

\subsection{Frobenius algebras}

Throughout this paper, 
$\bk$ will be a commutative $\bQ$-algebra.
An {\em invariant metric} on a commutative algebra  
$(H, \wedge)$ over $\bk$ is a non-degenerate bilinear map 
$\eta: H \times H \to \bk$ such that 
\begin{eqnarray}
&& \eta (X, Y) = \eta (Y, X), \label{eqn:symmetric} \\
&& \eta(X \wedge Y, Z) = \eta(X, Y \wedge Z), \label{eqn:Frobenius1}
\end{eqnarray}
for $X, Y, Z \in H$.
The triple $(H, \wedge, \eta)$ is 
called a {\em Frobenius algebra}.
Suppose that $H$ is free as $\bk$-module,
and fix a basis $\{e_a\}$ of $H$.
Then there are elements $\phi_{ab}^c \in \bk$,
such that 
$$e_a \wedge e_b = \phi_{ab}^c e_c.$$
Let $\eta_{ab} = \eta(e_a, e_b)$ and
$\phi_{abc} = \eta(e_a \wedge e_b, e_c) = \phi_{ab}^p\eta_{pc}$.
From (\ref{eqn:symmetric}) and (\ref{eqn:Frobenius1}),
one easily sees that $\phi_{abc}$ is symmetric in all three indices.
Denote by $(\eta^{ab})$ the inverse matrix of $(\eta_{ab})$,
then $\phi_{ab}^c = \phi_{abp}\eta^{pc}$.
The associativity of the multiplication is equivalent to the following
system of equations
\begin{eqnarray} \label{eqn:WDVV1}
\phi_{abp}\eta^{pq}\phi_{qcd} = \phi_{bcp}\eta^{pq}\phi_{aqd}.
\end{eqnarray}
When $(H, \wedge)$ has an identity $1$, 
only $\phi$ is needed.
In fact, 
one can take $e_0 = 1$.
Then from (\ref{eqn:Frobenius1}),
one gets $\eta_{ab} = \phi_{0ab}$.
To summarize,
the structure of a Frobenius algebra with identity is determined by the 
symmetric $3$-tensor $\phi$
such that $\eta_{ab} = \phi_{0ab}$.
One can easily generalize this discussion to the graded case.

\begin{example}
Let $M$ be an oriented connected closed $n$-dimensional smooth manifold,
then the de Rham cohomology $H^*(M)$ with the wedge product $\wedge$
is a graded commutative algebra over $\bR$.
Let $\int_M: H^*(M) \to \bR$ be defined by 
integrations of the top degree components.
Then Poincar\'{e} duality implies that
$\int_M$ induces a Frobenius algebra structure on $(H^*(M), \wedge)$. 
\end{example}

\subsection{WDVV equations and formal Frobenius manifold structures}
Denote by $\{x^a\}$ the linear coordinates in the basis $\{e_a\}$.
Consider a family of commutative 
associative multiplications $\wedge_x$ on $H$, 
one for each $x \in H$,
 such that 
$$\eta(X \wedge_x Y, Z) = \eta(X, Y \wedge_x Z),$$
for all $X, Y, Z, x \in H$.
We then have a family of $3$-tensors $\phi_{abc}(x)$.
Such families with the property that
$$\frac{\partial}{\partial x^d} \phi_{abc} 
= \frac{\partial}{\partial x^c} \phi_{abd}$$
are particularly interesting to physicists
(see e.g. Dijkgraaf-Verlinde-Verlinde \cite{Dij-Ver-Ver}).
Given such a family,
one can find a function $\Phi: H \to \bk$,
such that
$$\phi_{abc} 
= \frac{\partial^3\Phi}{\partial x^a \partial x^b\partial x^c}.$$
By (\ref{eqn:WDVV1}), 
$\Phi$ satisfies the Witten-Dijkgraaf-E. Verlinde-H. Verlinde
(WDVV) equations:
\begin{eqnarray} \label{eqn:WDVV}
\frac{\partial^3 \Phi}{\partial x^a \partial x^b \partial x^p} \eta^{pq}
\frac{\partial^3 \Phi}{\partial x^q \partial x^c \partial x^d} =
\frac{\partial^3 \Phi}{\partial x^b \partial x^c \partial x^p} \eta^{pq}
\frac{\partial^3 \Phi}{\partial x^a \partial x^q \partial x^d}.
\end{eqnarray}
The function $\Phi$ is called the {\em potential function} for the family.
A {\em formal Frobenius manifold structure} on
$(H, \wedge, \eta)$ is a formal power series $\Phi$ 
which satisfies the WDVV equations (Manin \cite{Man1}).
We refer to $(H, \wedge)$
as the initial data for the WDVV equations.
If $(H, \wedge)$ has an identify $1 = e_0$ which is also an identity for
all of $\wedge_x$,
then we have
\begin{eqnarray} \label{eqn:metric}
\eta_{ab} 
= \frac{\partial^3\Phi}{\partial x^0\partial x^a \partial x^b}.
\end{eqnarray}
If a formal Frobenius manifold structure $\Phi$ satisfies (\ref{eqn:metric}),
it is called a structure of formal Frobenius manifold with identity.
Again,
it is straightforward to generalize the above discussions to the graded case.

\section{Some notions from rational homotopy theory}

\subsection{Quasi-isomorphisms}  

A {\em quasi-isomorphism} between two 
DGA's $\cA$ and $\cB$ is a 
series of DGAs $\cA_0, \cdots, \cA_n$, 
and DGA-homomorphisms either $f_i: \cA_i \rightarrow \cA_{i+1}$, 
or $f_i: \cA_{i+1} \rightarrow \cA_i$ for $0 \leq i \leq n-1$,
such that $\cA_0 = \cA$, $\cA_n = \cB$ and each $f_i$ induces 
isomorphism on cohomology.
A DGA $\cA$ is called {\em formal} 
if it is quasi-isomorphic to its cohomology algebra
(regarded as a DGA with zero differential).
Every simply connected DGA has a minimal model and
quasi-isomorphic DGAs have the same minimal model
(see e.g. Giffiths-Morgan \cite{Gri-Mor}).

\subsection{Formal power series connections}
\label{sec:formalpower}

Given a  {\em differential graded Lie algebra} (DGLA) 
$(\cL, [\cdot, \cdot], \fd)$,
fix a decomposition $\cL = \cH \oplus \fd M \oplus M$,
such that $\cH \subset \Ker \fd$, 
the natural map $\cH \to H(\cL, \fd)$ is
an isomorphism and $\fd|_M$ is injective.
Such a decomposition is called a {\em cohomological decomposition}.
Fix a homogeneous basis $\{\alpha_j \in \cH\}$,
denote by $\{X^j\}$ the dual basis.
Each $X^j$ is given the degree $- |\alpha_j| + 1$.
The dual space of $\cH$ with such a grading is denoted by $s^{-1}\cH^t$.
Let $\overline{S}(s^{-1}\cH^t) = \sum S^n(s^{-1}\cH^t)$.
Modifying the method in Chen \cite{Che},
Hain \cite{Hai} inductively constructed a differential
$b$ on $\overline{S}(s^{-1}\cH^t)$ and a 
{\em formal power series connection} 
of the form
$$\omega = \sum \alpha_i X^i + \cdots 
+\sum \alpha_{i_1\cdots i_n} X^{i_1}\odot \cdots \odot X^{i_n} + \cdots,$$
where $\alpha_{i_1\cdots i_n} \in \cL$ 
has degree $1-|X^{i_1}| - \cdots - |X^{i_n}|$,
such that
$$b \omega + \fd \omega + \frac{1}{2} [\omega, \omega] = 0.$$
Here we have naturally extended $b$, $\fd$ and 
$[\cdot, \cdot]$ on $\cL \otimes \overline{S}(s^{-1}\cH^t)$.

\section{Formal Frobenius manifold structures from DGBV algebras}
\label{sec:DGBV}

\subsection{DGBV algebras}
\label{subsec:DGBV}

A {\em Gerstenhaber algebra} (G-algebra) is 
a graded commutative algebra
$(\cA, \wedge)$ over $\bk$ together with a bilinear map
$[\cdot \bullet \cdot]: \cA \otimes \cA \to \cA$ of degree $-1$,
such that
\begin{eqnarray*}
&& [a \bullet b] = - (-1)^{(|a|-1)(|b|-1)} [b \bullet a], \\
&& [a \bullet [b \bullet c]] = [[a \bullet b] \bullet c]
	+ (-1)^{(|a|-1|)(|b|-1)}[b \bullet [a \bullet c]], \\
&& [a \bullet (b \wedge c)] = [a \bullet b] \wedge c
	+ (-1)^{(|a|-1|)|b|} b \wedge [a \bullet c],
\end{eqnarray*}
for homogeneous $a, b, c \in \cA$.
For any linear operator $\Delta: \cA \to \cA$ of degree $-1$,
define
$$[a \bullet b]_{\Delta}= (-1)^{|a|} (\Delta (a \wedge b) 
- (\Delta a) \wedge b - (-1)^{|a|} a \wedge \Delta b),$$
for homogeneous elements $a, b \in \cA$.
If $\Delta^2 = 0$ and
\begin{eqnarray*}
[a \bullet (b \wedge c)]_{\Delta} 
	= [a \bullet b]_{\Delta} \wedge c
	+ (-1)^{(|a|-1)|b|} b \wedge [a \bullet c]_{\Delta},
\end{eqnarray*}
for all homogeneous $a, b, c \in \cA$,
then $(\cA, \wedge, \Delta, [\cdot \bullet \cdot]_{\Delta})$
is a {\em Gerstenhaber-Batalin-Vilkovisky (GBV) algebra}. 
When there is no confusion, 
we will simply write $[\cdot \bullet \cdot]$ 
for $[\cdot \bullet \cdot]_{\Delta}$.
It can be checked that a GBV algebra is a $G$-algebra
(cf. Koszul \cite{Kos}, Getzler \cite{Get} and Manin \cite{Man2}).
A  {\em DGBV algebra} is a GBV algebra with a differential $\delta$
with respect to $\wedge$, 
such that $\delta\Delta + \Delta \delta = 0$
(hence $\delta$ is also a differential of 
$[\cdot \bullet \cdot]_{\Delta}$).

\subsection{Nice integrals}
A $\bk$-linear functional $\int: {\cal A} \rightarrow \bk$ 
on a DGBV algebra $(\cA, \wedge, \delta, \Delta, [\cdot \bullet \cdot])$
is called {\em an integral} if
for all homogeneous $a, b \in {\cal A}$,
\begin{eqnarray}
\int (\delta a) \wedge b & = & 
	(-1)^{|a|+1} \int a \wedge \delta b, \label{int1} \\
\int (\Delta a ) \wedge b & = &
	(-1)^{|a|} \int a \wedge \Delta b. \label{int2}
\end{eqnarray} 
Given an integral,
$\eta (\alpha, \beta) = \int \alpha \wedge \beta$ defines 
a graded symmetric bilinear form $\eta$ on $H({\cal A}, \delta)$.
If $\eta$ is non-degenerate, 
we say that the integral is {\em nice}.
Therefore,
if ${\cal A}$ has a nice integral, 
then $(H(\cA, \delta), \wedge, \eta)$ is a Frobenius algebra.

\subsection{Formal Frobenius manifold structure on the cohomology of 
a DGBV algebra}
\label{sec:DGBV2Frobenius}

Let $({\cal A}, \wedge, \delta, \Delta, [\cdot \bullet \cdot])$
be a DGBV algebra satisfying the following conditions:
\begin{itemize}
\item[(i)] The cohomology algebra 
$H({\cal A}, \delta) = \Ker \delta /\Img \delta$ 
is free and of finite rank as a $\bk$-module.
\item[(ii)] There is a nice integral on ${\cal A}$.
\item[(iii)] The inclusions $(Ker \Delta, \delta) \hookrightarrow 
({\cal A}, \delta)$ and $(Ker \delta, \Delta) \hookrightarrow 
({\cal A}, \Delta)$  induce isomorphisms of cohomology.
\end{itemize}
Then there is a canonical construction of a
formal Frobenius manifold  structure with identity on $H({\cal A}, \delta)$.

This construction was implicit in 
Bershadsky-Cecotti-Ooguri-Vafa \cite{Ber-Cec-Oog-Vaf}
in the special case of the extended deformation theory of Calabi-Yau manifolds
based on the work of Tian \cite{Tia} and Todorov \cite{Tod}.
It was mathematically formulated by  Barannikov and Kontsevich \cite{Bar-Kon}.
The details for general DGBV algebras can be found in Manin \cite{Man2}.
Here we give a description in terms of Chen's construction of formal power 
series connections.
Since $(s\cA, \delta, [\cdot\bullet\cdot])$ is a DGLA,
by \S \ref{sec:formalpower}
one gets a universal formal power series connection 
$\Gamma$ and a derivation $\partial$ on 
$K = \overline{S}(H(\cA, \delta)^t)$,
such that 
$$\partial \Gamma + \delta \Gamma + \frac{1}{2}[\Gamma\bullet\Gamma] = 0.$$
However,
the condition (iii) above implies that $\partial = 0$. 
Such a method was first used in Tian \cite{Tia} and Todorov \cite{Tod}.
Therefore, $\Gamma$ satisfies
\begin{eqnarray}
\delta \Gamma + \frac{1}{2}[\Gamma \bullet \Gamma]= 0, 
	\label{eqn:MC} 
\end{eqnarray}
Furthermore,
one can take $\Gamma = \Gamma_1 + \Delta B$,
where $\Gamma_1 = \sum \alpha_j X^j$.
Set $\delta_{\Gamma} = \delta + [\Gamma\bullet \cdot]$,
then $\delta_{\Gamma}$ is a derivation of $\cA_K = (\cA \otimes K, \wedge)$.
Now
$$\delta_{\Gamma}^2 
= [(\delta \Gamma + \frac{1}{2}[\Gamma\bullet\Gamma])\bullet \cdot] = 0.$$
It can be proved that 
$H(\cA_K, \delta_{\Gamma}) \cong H(\cA, \delta) \otimes K$,
hence the multiplication in $H(\cA_K, \delta_{\Gamma})$ 
induces a deformation of the multiplication in $H(\cA, \delta)$.
Given $X \in H(\cA, \delta)$, 
the contraction with $X$ from the right induces a right
derivation of degree $|X|$ 
on $K$ and hence also on 
$\cA_K$.
For $\alpha \in \cA_K$,
denote by $X\alpha$ the contraction by $X$ of $\alpha$. 
Now applying contraction by $X$ on both sides of (\ref{eqn:MC}),
one gets
$$\delta (X\Gamma) + [\Gamma \bullet (X\Gamma)] = 0.$$
So $X\Gamma$ represents a class in $H(\cA_K, \delta_{\Gamma})$.
Following Bershadsky {\em et al} \cite{Ber-Cec-Oog-Vaf},
set
\begin{eqnarray} \label{eqn:KSaction}
\Phi = \int \frac{1}{6}\Gamma^3 - \frac{1}{2}\delta B \Delta B.
\end{eqnarray}
A calculation in Barannikov-Kontsevich \cite{Bar-Kon}
and Manin \cite{Man2} shows that 
$$X^3 \Phi = \int (X\Gamma) \wedge (X\Gamma) \wedge (X\Gamma),$$
i.e. $\Phi$ is the potential function of the deformation.

\subsection{Gauge invariance}
We recall some results proved in Cao-Zhou \cite{Cao-Zho5}.
Let $\fm$ be the maximal ideal of $K$.
Consider the group
\begin{eqnarray*}
\cG = \exp (s\cA \otimes \fm)^e,
\end{eqnarray*}
with the multiplication $e^A e^B = e^C$  defined
by the Campbell-Baker-Hausdorff formula.
Assume $(\cA, \wedge, \delta, \Delta, [\cdot\bullet\cdot])$ is a DGBV algebra
which satisfies the conditions (i)-(iii) in \S \ref{sec:DGBV2Frobenius}.
Then it is shown in \cite{Cao-Zho5} that given any two universal 
normalized solutions $\Gamma$ and $\overline{\Gamma}$,
there exists an odd element $A \in (\Img \Delta \otimes \fm)$, 
such that $e^A \cdot \Gamma = \overline{\Gamma}$.
Furthermore,
the potential function $\Phi$ is gauge invariant:
$\Phi(e^A \cdot \Gamma) = \Phi(\Gamma)$.

\subsection{Invariance under quasi-isomorphisms}

A homomorphism between two DGBV algebras with nice integrals over $\bk$,
 $(\cA_i, \wedge_i, \delta_i, \Delta_i, [\cdot\bullet\cdot]_i, \int_i)$,
$i = 1, 2$, 
is a homomorphism of graded algebras $f: \cA_1 \to \cA_2$
such that $f \delta_1 = \delta_2 f$, $f \Delta_1 = \Delta_2 f$,
and $\int_2 f(\alpha) = \int_1 \alpha$ for all $\alpha \in \cA$.
It is a quasi-isomorphism if $f$ induced isomorphisms 
$H(\cA_1, \delta_1) \to H(\cA_2, \delta_2)$
and $H(\cA_2, \Delta_2) \to H(\cA_2, \Delta_2)$.
We prove in \cite{Cao-Zho5} the following result:
If there is a quasi-isomorphism 
between two DGBV algebras with nice integrals
which satisfy the conditions (i)-(iii) in \S \ref{sec:DGBV2Frobenius},
then the formal Frobenius manifolds constructed from them 
as in \S \ref{sec:DGBV2Frobenius}
can be naturally identified with each other.

\section{DGBV algebras from symplectic and complex geometries}
\label{sec:DGBVfromPC}

\subsection{GBV algebra structure on the space of polyvector fields}
Denote by $\Omega^{-*}(M) = \Gamma(M, \Lambda^*TM)$ the space of polyvector 
fields on $M$.
The Lie bracket $[\cdot, \cdot]$ 
of vector fields can be extended to the {\em Schouten-Nijenhuis bracket}
$[\cdot, \cdot]_S: \Omega^{-*}(M) \times \Omega^{-*}(M) \to \Omega^{-*}(M)$
of degree $-1$ given locally by 
\begin{eqnarray*}
&& [X_1 \wedge \cdots \wedge X_p, Y_1 \wedge \cdots \wedge Y_q]_S \\
& = & \sum_{i = 1}^p\sum_{j=1}^q
(-1)^{i+j+p+1} [X_i, Y_j] \wedge 
X_1 \wedge \cdots \wedge \widehat{X}_i \wedge \cdots \wedge X_p
\wedge Y_1 \wedge \cdots \wedge \widehat{Y}_j \wedge \cdots \wedge Y_q, 
\end{eqnarray*}
for local vector fields $X_1, \cdots, X_p, Y_1, \cdots, Y_q$.
Then $(\Omega^{-*}(M), \wedge, [\cdot, \cdot]_S)$ becomes a G-algebra.
There are two methods to make it a GBV algebra.
The first method due to Koszul \cite{Kos} 
uses the generalized divergence operator 
for any torsion free connection, while
the second method due to Witten \cite{Wit} uses any volume form in a way
similar to earlier construction in Tian \cite{Tia} and Todorov \cite{Tod}
(see \S \ref{sec:CY}).

\subsection{DGBV algebras from Poisson geometry}
\label{sec:Poisson}
A {\em Poisson structure} on a manifold $M$ is an element 
$w \in \Omega^{-2}(M)$ such that $[w, w]_S = 0$.
On a Poisson manifold $(M, w)$,
define $\sigma: \Omega^{-*}(M) \to \Omega^{-(*+1)}(M)$ by
$\sigma(P) = [w, P]$.
Then $(\Omega^{-*}(M), \wedge, \sigma, [\cdot, \cdot]_S)$ is a differential
G-algebra.
If a  Poisson manifold $(M, w)$
is {\em regular}, i.e.,  $w$ has constant rank, then there is a (not unique)
torsion free connection for which the Poisson structure $w$ is parallel.
If $\Delta $ is the generalized divergence operator for such a connection, 
then $[\Delta, \sigma] = 0$ and
$(\Omega^{-*}(M), \wedge, \sigma, \Delta, [\cdot, \cdot]_S)$
is a DGBV algebra.

There is another DGBV algebra associated with a Poisson manifold $(M, w)$.
Koszul defined an operator 
$\Delta = [i(w), d]: \Omega^*(M) \to \Omega^{*-1}(M)$,
where $i(w)$ denotes the contraction by $w$.
He proved $[d, \Delta] = \Delta^2 = 0$ and 
$[\cdot\bullet\cdot]_{\Delta}$ satisfies
$$[\alpha\bullet(\beta \wedge\gamma)]_{\Delta}
=[[\alpha\bullet\beta]_{\Delta} \wedge \gamma
+ (-1)^{(|\alpha|-1)|\beta|}\beta \wedge [\alpha \bullet \gamma]_{\Delta},
$$
hence $(\Omega^*(M), \wedge, d, \Delta, [\cdot\bullet\cdot]_{\Delta})$
is a DGBV algebra.

\subsection{DGBV algebras from Symplectic manifolds}
\label{sec:symplectic}

A {\em symplectic manifold} $M$ is automatically a Poisson manifold,
hence there are associated DGBV algebras from it.
Indeed,
the symplectic form $\omega$ induces isomorphisms 
$\sharp: \Omega^*(M) \to \Omega^{-*}(M)$
and $\flat: \Omega^{-*}(M) \to \Omega^*(M)$.
Using the Darboux theorem,
it can be easily seen that the bivector filed $w = \omega^{\sharp}$ is a
Poisson structure.
Furthermore,
it can be seen that $d = [\omega\bullet\cdot]_{\Delta}$.
The symplectic form induces an isomorphism of DGBV algebras
$(\Omega^{-*}(M), \wedge, \sigma, \Delta, [\cdot, \cdot]_S)
\cong (\Omega^*(M), \wedge, d, \Delta, [\cdot\bullet\cdot]_{\Delta})$.

\subsection{Differential G-algebras from complex manifolds}
\label{sec:complex}

Given a complex $n$-manifold $M$,
let 
\begin{eqnarray*}
\Omega^{-*, -*}(M) = \sum_{p, q\geq 0} \Omega^{-p, -q}(M) 
= \sum_{p, q\geq 0}  
\Gamma(M, \Lambda^pT^{1, 0}M \otimes \Lambda^qT^{0,1}M).
\end{eqnarray*}
Similarly define $\Omega^{-*, *}(M)$ and $\Omega^{*, *}(M)$.
Since $\Omega^{-*, -*}(M) = \Omega^{-*}(M) \otimes \bC$,
it is a complex G-algebra.
From $[\Omega^{-1, 0}(M), \Omega^{-1, 0}(M)] \subset \Omega^{-1, 0}(M)$,
one sees that $\Omega^{-*, 0}(M)$ is G-subalgebra of
$\Omega^{-*, -*}(M)$.
Using local coordinates,
it is easy to see that the Schouten-Nijenhuis bracket on
$\Omega^{-*, 0}(M)$ can be extended to a G-algebra 
structure on $\Omega^{-*, *}(M)$.
Since $\Omega^{-*, *}(M)$ is the space of $(0, *)$-form 
with values in the holomorphic bundle $\Lambda^*T^{1, 0}M$, 
there is an operator $\bar{\partial}$ acting on it.
The tuple $(\Omega^{-*, *}(M), \wedge, [\cdot, \cdot]_S, \bar{\partial})$ is 
a differential G-algebra.

\subsection{DGBV algebras from Calabi-Yau manifolds}
\label{sec:CY}

Recall that a K\"{a}hler manifold is called a {\em Calabi-Yau manifold}
if it admits a parallel holomorphic volume form.
The terminology comes from Yau's solution to Calabi's conjecture 
\cite{Yau2}.
Such manifolds are very important in string theory.
Given a Calabi-Yau $n$-fold $M$,
the holomorphic volume form $\Omega$
defines an isomorphism $\Omega^{-*, *}(M) \to \Omega^{n-*, *}(M)$.
Let $\Delta: \Omega^{-*, *}(M) \to \Omega^{-(*-1), *}(M)$
be the conjugation of $\partial$ by this isomorphism.
A formula in Tian \cite{Tia} shows that 
$[\cdot\bullet\cdot]_{\Delta} = [\cdot, \cdot]_S$,
hence $\Delta$ gives a GBV algebra structure on $\Omega^{-*, *}(M)$.
From $\bar{\partial}\Omega = 0$,
one sees that $\bar{\partial} \Delta + \Delta \bar{\partial} = 0$.
Hence $(\Omega^{-*, *}(M), \wedge, \bar{\partial}, \Delta, [\cdot, \cdot])$
is a DGBV algebra.

\section{Deformations of complex and Poisson structures}

\subsection{Deformation of complex structures}

Given a complex manifold $M$,
the small deformations of the complex structure can be studied via 
the Maurer-Cartan equation
\begin{eqnarray} \label{eqn:complexMC}
\bar{\partial} \omega + \frac{1}{2} [\omega, \omega] = 0,
\end{eqnarray}
where $\omega \in \Omega^{-1, 1}(M)$.
Fix a Hermitian metric on $M$,
let $\bar{\partial}^*$ be the formal adjoint of $\bar{\partial}$
and $\square_{\bar{\partial}} 
= \bar{\partial}\bar{\partial}^* + \bar{\partial}^*\bar{\partial}$.
As a consequence of the Hodge decomposition 
$$\Omega^{-*, *}(M) 
= \cH^{-*, *}_{\bar{\partial}}(M) \oplus \Img \bar{\partial}
\oplus \Img \bar{\partial}^*,$$
we have $H^{-*, *}(M) \cong \cH^{-*, *}_{\bar{\partial}}(M)$.
Take a basis $\{\alpha_j\}$ of $\cH^{-1, 1}_{\bar{\partial}}$,
and let $\{t^j\}$ be coordinates in this basis.
Let $\omega(t) = \omega_1(t) + \cdots + \omega_n(t) + \cdots$ 
be a formal power series,
such that $\omega_1(t) = \sum_j \alpha_j t^j$,
and $\omega_n$ is homogeneous of degree $n$ in $t^j$'s.
Then (\ref{eqn:complexMC}) is equivalent to a sequence of equations
$$\bar{\partial}\omega_n = - \frac{1}{2} \sum_{p+q=n} [\omega_p, \omega_q].$$
Suppose that $\omega_1, \cdots, \omega_n$ have been defined,
a calculation shows that 
$$\bar{\partial} (- \frac{1}{2} \sum_{p+q=n+1} [\omega_p, \omega_q]) = 0.$$
Let $Q = G_{\bar{\partial}}\bar{\partial}^*$,
where $G_{\bar{\partial}}$ is the Green's operator 
of $\square_{\bar{\partial}}$.
Take 
$$\omega_{n+1} = - \frac{1}{2} Q\sum_{p+q=n+1} [\omega_p, \omega_q].$$
This is equivalent to
\begin{eqnarray} \label{eqn:Kuranishi}
\omega + \frac{1}{2}Q[\omega, \omega] = \omega_1.
\end{eqnarray}
Following Kodaira-Spencer \cite{Kod-Spe},
one can show that for small $||t||$,
$\omega(t)$ is convergent.
Taking $\bar{\partial}$ on both sides of (\ref{eqn:Kuranishi}),
one gets
$$\bar{\partial}\omega + \frac{1}{2} [\omega, \omega] 
= -\frac{1}{2} [\omega, \omega]^H,$$
where $[\omega, \omega]^H$ is the harmonic part of $[\omega, \omega]$.
When $H^{-1, 2}(M) = 0$,
a neighborhood of $0$ in $H^{-1, 1}(M)$ then parameterizes 
the small deformations.
This is the construction by Kodaira-Spencer \cite{Kod-Spe} 
of a complete family when $H^{-1,2}(M) = 0$.
When  $H^{-1, 2}(M) \neq 0$,
one gets a map from $ U \subset \cH^{-1, 1}_{\bar{\partial}}(M)$ 
to $\cH^{-1, 2}_{\bar{\partial}}(M)$ given by 
$t = (t^j)\mapsto [\omega(t), \omega(t)]^H$,
its zero set gives solutions $\omega$ to the Maurer-Cartan equation.
This is due to Kuranishi (see e.g. \cite{Kur}).

The extended deformation problem in this case is to 
find solutions to the Maurer-Cartan equation with 
$\omega \in \Omega^{-*, *}(M)$.
The Hodge decomposition gives the cohomological splitting, and
the above iterative method of $\omega$ corresponds to Chen's construction
of formal power series connection modified by
Hain \cite{Hai} for DGLA's.
The formal power series 
$-\frac{1}{2}[\omega(t), \omega(t)]^H$ corresponds to the differential $b$
(see \S \ref{sec:formalpower}), 
where $\{t^j\}$ are the coordinates of $\cH^{-*, *}(M)$ 
in a homogeneous basis.

\subsection{Formal Frobenius manifolds from Calabi-Yau manifolds}
Motivated by a result of Bogomolov \cite{Bog},
Tian \cite{Tia} (see also Todorov \cite{Tod}) 
introduced an ingenious method to prove that the deformations
of complex structures on Calabi-Yau manifolds are unobstructed.
This includes the introduction of the operator $\Delta$ in \S 
\ref{sec:CY}, showing $[\cdot\bullet\cdot]_{\Delta} = [\cdot, \cdot]_S$,
using Hodge theory to find power series solution $\omega(t)$ such that
$\omega(t) = \omega_1(t) + \Delta B(t)$,
and then showing that $[\omega(t), \omega(t)]^H = 0$.
This method has been used by 
Bershadsky-Cecotti-Ooguri-Vafa \cite{Ber-Cec-Oog-Vaf}
in what they called Kodaira-Spencer theory of gravity.
Fixing a holomorphic volume form $\Omega$ on a Calabi-Yau $n$-manifold $M$,
define $\int: \Omega^{-*, *}(M) \to \bC$ by
$$\int \alpha = \int_M \alpha^{\flat} \wedge \Omega.$$
It is easy to see that
$\int$ is a nice integral for 
$(\Omega^{-*,*}(M), \wedge, \bar{\partial}, \Delta, [\cdot, \cdot]_{\Delta})$.
Bershadsky {\em et al.} found a Lagrangian invariant under the action of 
diffeomorphisms.
They showed that its critical points correspond to 
the complex structures nearby. 
They discussed the quantization via BV formalism
in which the extended deformation problem suggested by Witten arises naturally. 
They also argued how the values of the Lagrangian 
at the critical points give the potential
for the deformed multiplicative structure.
Barannikov and Kontsevich \cite{Bar-Kon} 
formulated such results in terms of Frobenius manifolds
introduced by Dubrovin \cite{Dub}.
They also remarked that such construction 
can be carried out for DGBV algebras satisfying 
conditions  in \S \ref{sec:DGBV2Frobenius}.
The details of this construction of formal Frobenius manifold
can be found in Manin \cite{Man2}.
Anyway,
the moral is that the extended deformation problem 
leads to the DGBV algebra
$(\Omega^{-*, *}(M), \wedge, \bar{\partial}, \delta, [\cdot, \cdot]_S)$,
which further leads to a canonical construction of formal Frobenius manifold 
structure
on $H^{-*, *}(M)$.

\subsection{Deformations of Poisson structures}
Given a Poisson structure $w_0$ on $M$,
any Poisson structure close to $w_0$ can be written as $w = w_0 + \gamma$.
Let $\sigma = [w_0, \cdot]_S: \Omega^{-*}(M) \to \Omega^{-(*+1)}(M)$,
then $[w, w]_S = 0$ can be rewritten as 
\begin{eqnarray} \label{eqn:PoissonMC1}
\sigma\gamma + \frac{1}{2}[\gamma, \gamma]_s = 0.
\end{eqnarray}
Now assume that $w_0$ comes from a symplectic form $\omega_0$ on $M$.
The deformation theory of symplectic structure is ``trivial".
It is well-known that any closed $2$-form close to a symplectic $2$-form  
in $C^0$-norm is also symplectic.
However,
the deformation of the Poisson structure is nonlinear:
equation (\ref{eqn:PoissonMC1}) corresponds to 
\begin{eqnarray} \label{eqn:PoissonMC2} 
d \omega + \frac{1}{2} [\omega \bullet \omega] = 0.
\end{eqnarray}
Given a Riemannian metric on $M$,
then one can consider the formal adjoint $d^*$ of $d$ 
and the Laplacian operator $\square_d: \Omega^*(M) \to \Omega^*(M)$.
The power series method can be carried out similar to the complex case.
Also one can consider 
the extended deformation problem of finding solutions of 
(\ref{eqn:PoissonMC2}) in $\Omega^*(M)$.

\subsection{K\"{a}hler gravity}
The mirror analogue of the Kodaira-Spencer theory of gravity 
is the theory of K\"{a}hler gravity studied 
in Bershadsky-Sadov \cite{Ber-Sad}.
Let $M$ be a closed K\"{a}hler manifold with K\"{a}hler form $\omega_0$.
By Stokes theorem and Poincar\'{e} duality, 
the integral of top degree forms on $M$
is a nice integral of the DGBV algebra 
$(\Omega^*(M), \wedge, d, \Delta, [\cdot\bullet\cdot])$.
K\"{a}hler identities shows that $\Delta = (d^c)^*$.
By imposing the gauge fixing condition $\Gamma \in \Ker \Delta$,
Bershadsky and Sadov wrote a Lagrangian similar 
to that used in Kodaira-Spencer gravity
and showed that
the Euler-Lagrangian equation of it on $\Ker \Delta$ is
$$d \Gamma + \frac{1}{2}\Delta(\Gamma \wedge \Gamma) = 0,$$
or equivalently,
$$d \Gamma + \frac{1}{2} [\Gamma\bullet\Gamma] = 0.$$
According to \S \ref{sec:symplectic}, 
it describes the deformations of the Poisson structure
corresponding to $\omega_0$.
Given a solution $K_0$ with $(d^c)^*K_0 = 0$,
Bershadsky and Sadov pointed out the operator 
$D = d + [(d^c)^*, K_0]$ squares to zero.
Notice that $D = d + [K_0 \bullet\cdot] = [(\omega_0 + K_0) \bullet \cdot]$,
so deforming $d$ to $D$ corresponds to deforming $[w_0, \cdot]_S$ to 
$[w, \cdot]_S$ on $\Omega^{-*}(M)$.

Similarly,
consideration of BV quantization leads to the extended deformation problem,
which is described by 
the DGBV algebra $(\Omega^*(M), \wedge, d, \Delta, [\cdot\bullet\cdot])$.
Again, Hodge theory of K\"{a}hler manifold can be used to show that this DGBV 
algebra
satisfies all the conditions in \S \ref{sec:DGBV2Frobenius}.
Hence there is a canonical construction of 
formal Frobenius manifold structure on the de Rham cohomology $H^*(M)$ for any 
K\"{a}hler manifold.
This was done in Cao-Zhou \cite{Cao-Zho2}.
Notice that we can complexify this DGBV algebra.
Then we have $d = \bar{\partial} + \bar{\partial}$,
$(d^c)^* = \sqrt{-1} \bar{\partial}^* - \sqrt{-1}\partial^*$.
In Cao-Zhou \cite{Cao-Zho1},
we proved that
$(\Omega^{*, *}(M), \wedge, \bar{\partial}, -\sqrt{-1}\partial^*,
[\cdot\bullet\cdot]_{-\sqrt{-1}\partial^*})$ is a DGBV algebra
which satisfies the conditions in \S \ref{sec:DGBV2Frobenius}.
Hence there is a natural construction of the Frobenius manifold structure on 
the Dolbeault cohomology $H^{*, *}(M)$.
Furthermore,
we showed in \cite{Cao-Zho2} that 
the formal Frobenius manifold structures on $H^*(M, \bC)$ and $H^{*, *}(M)$
given above can be identified with each other.
This generalizes of the well-known isomorphism 
$H^*(M, \bC) \cong H^{*, *}(M)$ (as complex vector spaces) 
to the highly nonlinear formal Frobenius manifold structures
on these spaces.

The above results have been generalized in various directions.
When $M$ is a hyperk\"{a}hler manifold,
we showed in \cite{Cao-Zho3} 
there are many different ways to construct DGBV algebras
for which the construction in \S \ref{sec:DGBV2Frobenius} applies.
Also,
when a K\"{a}hler manifold admits a Hamiltonian action 
by holomorphic isometries,
we showed in \cite{Cao-Zho4}
that the Cartan model admits a structure of DGBV algebra 
for which the method of \S \ref{sec:DGBV2Frobenius} yields a formal Frobenius
manifold structure.

\section{Concluding remarks}

In the above we have discussed the relevance of rational homotopy theory
in an A-type theory, the theory of K\"{a}hler gravity,
and a B-type theory, the Kodaira-Spencer theory of gravity.
The strategy is to obtain DGBV algebra structures by considering extended 
deformation problems of the relevant geometric objects.
Ideas from rational homotopy theory can then be used to construct and 
identify formal Frobenius manifold structures on the the cohomology of
such algebras.
As a result,
we propose the following version of mirror symmetry conjecture:
For a Calabi-Yau manifold $M$ and a suitably defined mirror manifold 
$\widehat{M}$,
the DGBV algebras 
$$(\Omega^{*, *}(M), \wedge, d, (d^c)^*, [\cdot\bullet\cdot]_{(d^c)^*})$$
and
$$(\Omega^{-*, *}(\widehat{M}), \wedge, \bar{\partial}, 
\Delta, [\cdot, \cdot]_S)$$
are quasi-isomorphic,
hence by \cite{Cao-Zho5} the formal Frobenius manifold structures on 
their cohomology can be identified with each other.
The examples of mirror manifold constructed 
by Greene and Plesser \cite{Gre-Ple} 
involves the quotients of Calabi-Yau manifold by finite automorphism groups.
The extension of the results described above to the quotient case is in 
progress.
Finally we want to mention that the above version of mirror symmetry
suggests an application being developed of 
Quillen's closed model category theory to the study of mirror symmetry.

\end{document}